\documentclass[twoside]{amsart}
\usepackage{amssymb}
\usepackage{amscd}
\usepackage[french]{babel}

\newcommand{\ad}{\operatorname{ad}}

\newcommand{\rg}{\operatorname{rg}}

\newcommand{\card}{\operatorname{card}}
\begin{document}
\setcounter{page}{1}
\thispagestyle{empty}

\markboth{P. Tauvel, R.W.T. Yu}{P. Tauvel, R.W.T. Yu}


\begin{center}
{\bf {\Large SUR L'INDICE DE CERTAINES ALG\`{E}BRES DE LIE}}
\end{center}

\medskip

\begin{center}
{\bf \large par Patrice TAUVEL et Rupert W.T. YU}
\end{center}

\renewcommand{\thefootnote}{}
\footnote{{\em Mots-cl\'{e}s~:} Alg\`{e}bres de Lie semi-simples $-$
Sous-alg\`{e}bres sp\'{e}ciales $-$ Indice.

{\em Classification math.}~: 17B05 $-$ 17B20.}
\bigskip

{\small {\sc R\'{e}sum\'{e}.}~On donne une majoration de l'indice de
certaines alg\`{e}bres de Lie introduites dans [2] et [7]. On en
d\'{e}duit la preuve d'une conjecture de D. Panyushev formul\'{e}e
dans [7]. Nous formulons aussi une conjecture concernant l'indice de
ces alg\`{e}bres, et la prouvons dans des cas particuliers. Enfin,
nous donnons un r\'{e}sultat concernant l'indice des sous-alg\`{e}bres
paraboliques d'une alg\`{e}bre de Lie semi-simple.}

\medskip

{\small {\sc Abstract.}~We give an upper bound for the index of
certain Lie algebras, called of seaweed type, introduced in [2] and
[7]. We deduce from this a conjecture of D. Panyushev stated in
[7]. We conjecture that our upper bound gives in fact a formula for
the index of such Lie algebras, and we prove it in certain
cases. Finally, we obtain a result concerning the index of parabolic
subalgebras of a semisimple Lie algebra.}

\bigskip

{\bf 1.~Notations et rappels}

\bigskip

{\bf 1.1.}~Dans la suite, $\Bbbk$ est un corps commutatif
alg\'{e}briquement clos de caract\'{e}ristique nulle. Les alg\`{e}bres
de Lie consid\'{e}r\'{e}es sont d\'{e}finies et de dimension finie sur
$\Bbbk$. Si $V$ est un $\Bbbk$-espace vectoriel de dimension finie, on
le munit de la topologie de Zariski. On renvoie \`{a} [1], [3] et [8]
pour les concepts g\'{e}n\'{e}raux utilis\'{e}s.

\medskip

{\bf 1.2.}~Soient $\mathfrak{g}$ une alg\`{e}bre de Lie,
$\mathfrak{g}^{*}$ son dual. L'alg\`{e}bre $\mathfrak{g}$ op\`{e}re
dans $\mathfrak{g}^{*}$ au moyen de la repr\'{e}sentation
coadjointe. Ainsi, si $X,Y \in \mathfrak{g}$ et $f \in
\mathfrak{g}^{*}$, on a :
$$
(X.f)(Y) = f([Y,X]). 
$$

Avec les notations pr\'{e}c\'{e}dentes, on d\'{e}finit une forme
bilin\'{e}aire altern\'{e}e $\Phi_{f}$ sur $\mathfrak{g}$ en posant :
$$
\Phi_{f}(X,Y) = f([X,Y]).
$$
Le noyau de $\Phi_{f}$ est not\'{e} $\mathfrak{g}^{(f)}$. L'entier
$\chi (\mathfrak{g}) = \inf \{ \dim \mathfrak{g}^{(f)}; \, f \in
\mathfrak{g}^{*} \}$ est appel\'{e} l'indice de $\mathfrak{g}$. On dit
que $f$ est r\'{e}gulier si $\dim \mathfrak{g}^{(f)} = \chi
(\mathfrak{g})$ ; l'ensemble $\mathfrak{g}_{r}^{*}$ des
\'{e}l\'{e}ments r\'{e}guliers de $\mathfrak{g}^{*}$ est un ouvert non
vide de $\mathfrak{g}^{*}$.

\medskip

{\bf 1.3.}~Dans la suite, on suppose que $\mathfrak{g}$ est
semi-simple, de forme de Killing $\kappa$. On fixe une
sous-alg\`{e}bre de Cartan $\mathfrak{h}$ de $\mathfrak{g}$, et on
note $R$ le syst\`{e}me de racines du couple $(\mathfrak{g},
\mathfrak{h}$). Pour $\alpha \in R$, $\mathfrak{g}^{\alpha}$ est
l'espace radiciel associ\'{e} \`{a} $\alpha$. On d\'{e}signe par
$H_{\alpha}$ l'unique \'{e}l\'{e}ment de $[\mathfrak{g}^{\alpha},
\mathfrak{g}^{-\alpha}]$ tel que $\alpha (H_{\alpha}) = 2$. Si
$\lambda \in \mathfrak{h}^{*}$, on \'{e}crit $\langle \lambda ,
\alpha^{\vee}\rangle$ pour $\lambda (H_{\alpha})$. Si $P$ est une
partie de $R$, on pose :
$$
\mathfrak{g}^{P} = \operatornamewithlimits{\textstyle \sum}_{\alpha
\in P} \mathfrak{g}^{\alpha}.
$$

Soient $\Pi$ une base de $R$ et $R_{+}$ (resp. $R_{-}$) l'ensemble des
racines positives (resp. n\'{e}gatives) correspondant. On d\'{e}signe
par $\mathfrak{b}$ et $\mathfrak{b}_{-}$ les sous-alg\`{e}bres de
Borel de $\mathfrak{g}$ d\'{e}finies par :
$$
\mathfrak{b} = \mathfrak{h} \oplus \mathfrak{g}^{R_{+}} \ , \
\mathfrak{b}_{-} = \mathfrak{h} \oplus \mathfrak{g}^{R_{-}}.
$$

Pour toute partie $S$ de $\Pi$, on note $\mathbb{Z}S$
(resp. $\mathbb{N}S$) l'ensemble des combinaisons lin\'{e}aires \`{a}
coefficients entiers (resp. \`{a} coefficients entiers positifs) des
\'{e}l\'{e}ments de $S$, et on pose :
$$
R^{S} = R \cap \mathbb{Z}S \ , \ R_{+}^{S} = R_{+} \cap \mathbb{Z}S =
R \cap \mathbb{N}S \ , \ R_{-}^{S} = R^{S} \backslash R_{+}^{S}.
$$
L'ensemble $R^{S}$ est un syst\`{e}me de racines dans le sous-espace
vectoriel de $\mathfrak{h}^{*}$ qu'il engendre, $S$ est une base de
$R^{S}$, et $R_{+}^{S}$ (resp. $R_{-}^{S}$) est l'ensemble des racines
positives (resp. n\'{e}gatives) associ\'{e}. Si $S$ est une partie
connexe de $\Pi$, le syst\`{e}me $R^{S}$ est irr\'{e}ductible, et on
note $\varepsilon_{S}$ sa plus grande racine (relativement \`{a}
$S$). 

\medskip

{\bf 1.4.}~Soit $S$ une partie connexe de $\Pi$. Pour tout racine
$\alpha \in R_{+}^{S} \backslash \{ \varepsilon_{S}\}$, on sait que
$\langle \alpha , \varepsilon_{S}^{\vee}\rangle \in \{ 0,1\}$. Si
$R_{0}^{S}$ est l'ensemble des $\alpha \in R^{S}$ v\'{e}rifiant
$\langle \alpha , \varepsilon_{S}^{\vee}\rangle = 0$, on a $R_{0}^{S}
= R^{T}$, o\`{u} $T = S \cap R_{0}^{S}$. D'autre part, pour $\alpha
\in R_{+}\cap R_{0}^{S}$, on a $\alpha \pm \varepsilon_{S} \notin R$
([1], proposition 9, p. 149) donc, si $\alpha \ne \varepsilon_{S}$,
les racines $\alpha$ et $\varepsilon_{S}$ sont fortement orthogonales.

\medskip

{\bf 1.5.}~Rappelons une construction et quelques propri\'{e}t\'{e}s
d'un certain ensemble de racines deux \`{a} deux fortement
orthogonales ([5] et [6]).

Soit $S$ une partie de $\Pi$. Par r\'{e}currence sur le cardinal de
$S$, on d\'{e}finit un sous-ensemble $\mathcal{K}(S)$ de parties de
$S$ de la mani\`{e}re suivante :

a) $\mathcal{K}(\emptyset ) = \emptyset$.

b) Si $S_{1}, \dots , S_{r}$ sont les composantes connexes de $S$,
alors :
$$
\mathcal{K}(S) = \mathcal{K}(S_{1}) \cup \cdots \cup \mathcal{K}(S_{r}).
$$

c) Si $S$ est connexe, on a :
$$
\mathcal{K}(S) = \{ S\} \cup \mathcal{K}( \{ \alpha \in S; \, \langle
\alpha , \varepsilon_{S}^{\vee}\rangle = 0\} ).
$$

\medskip

{\bf 1.6.}~Il est clair que le cardinal de $\mathcal{K}(\Pi )$ ne
d\'{e}pend que de $\mathfrak{g}$, mais pas de $\mathfrak{h}$ et
$\Pi$. On note cet entier $k_{\mathfrak{g}}$. On donne dans le tableau
suivant la valeur de $k_{\mathfrak{g}}$ pour les diff\'{e}rents types
d'alg\`{e}bres de Lie simples. On a not\'{e} $[r]$ la partie
enti\`{e}re d'un nombre rationnel $r$.

\bigskip

\begin{center}
\begin{tabular}{|c|c|c|c|c|c|c|c|c|c|}\hline
\strut\rule[-1ex]{0ex}{4ex} & $A_{\ell}, \ell \geqslant 1$ &
$B_{\ell}, \ell \geqslant 2$ & $C_{\ell}, \ell \geqslant 3$ &
$D_{\ell} , \ell \geqslant 4$ & $E_{6}$ & $E_{7}$ & $E_{8}$ & $F_{4}$
& $G_{2}$ \\ \hline
$k_{\mathfrak{g}}$ \strut\rule[-2.8ex]{0ex}{7ex}& $\left[
\dfrac{\ell+1}{2}\right]$ & $\ell$ & $\ell$ & $2 \left[
\dfrac{\ell}{2}\right]$& $4$ & $7$ & $8$ & $4$ & $2$ \\ \hline
\end{tabular}
\end{center}

\bigskip

{\bf 1.7.~Lemme.}~{\em Soient $S$ une partie de $\Pi$ et $K,K' \in
\mathcal{K}(S)$. Alors : 

{\em (i)} $K$ est une partie connexe de $\Pi$.

{\em (ii)} On a ou $K \subset K'$, ou $K' \subset K$, ou $K$ et $K'$
sont des parties disjointes de $S$ telles que $\alpha + \beta \notin
R$ pour $\alpha \in R^{K}$ et $\beta \in R^{K'}$.

{\em (iii)} Si $K \ne K'$, $\varepsilon_{K}$ et $\varepsilon_{K'}$
sont fortement orthogonales. Par cons\'{e}quent, $\{ \varepsilon_{K};
\, K \in \mathcal{K}(S) \}$ est un ensemble de racines deux \`{a} deux
fortement orthogonales de $R$.}

\medskip

{\bf 1.8.}~Si $S$ est une partie de $\Pi$ et si $K \in
\mathcal{K}(S)$, on pose :  
$$
\Gamma^{K} = \{ \alpha \in R^{K}; \, \langle \alpha ,
\varepsilon_{K}^{\vee} \rangle > 0 \} \ , \ \Gamma_{0}^{K} =
\Gamma_{K} \backslash \{ \varepsilon_{K} \} \ , \ \mathfrak{a}_{K} =
\operatornamewithlimits{\textstyle \sum}_{\alpha \in \Gamma^{K}}
\mathfrak{g}^{\alpha} .
$$

\medskip

{\bf Lemme.}~{\em Soient $K,K' \in \mathcal{K}(S)$, $\alpha , \beta
\in \Gamma^{K}$, et $\gamma \in \Gamma^{K'}$.

{\em (i)} On a $\Gamma^{K} = R_{+}^{K} \backslash \{ \delta \in
R_{+}^{K}; \, \langle \delta , \varepsilon_{K}^{\vee}\rangle = 0\}$.

{\em (ii)} L'ensemble $R_{+}^{S}$ est la r\'{e}union disjointe des
$\Gamma^{K''}$ pour $K'' \in \mathcal{K}(S)$, et $\mathfrak{a}_{K}$
est une alg\`{e}bre de Heisenberg de centre
$\mathfrak{g}^{\varepsilon_{K}}$. 

{\em (iii)} Si $\alpha + \beta \in R$, alors $\alpha + \beta =
\varepsilon_{K}$. 

{\em (iv)} Si $\alpha + \gamma \in R$, alors ou $K \subset K'$ et
$\alpha + \gamma \in \Gamma^{K'}$, ou $K' \subset K$ et $\alpha +
\gamma \in \Gamma^{K}$.}

\medskip

{\bf 1.9.}~On d\'{e}signe par $E_{S}$ le sous-espace vectoriel de
$\mathfrak{h}^{*}$ engendr\'{e} par les $\varepsilon_{K}$, $K \in
\mathcal{K}(S)$. Ces racines \'{e}tant deux \`{a} deux fortement
orthogonales, la dimension de $E_{S}$ est \'{e}gale \`{a}
$\card (\mathcal{K}(S))$ (cardinal de $\mathcal{K}(S)$).

\bigskip

{\bf 2.~Sous-alg\`{e}bres sp\'{e}ciales}

\bigskip

{\bf 2.1.}~On d\'{e}signe toujours par $\mathfrak{g}$ une alg\`{e}bre
de Lie semi-simple, et on note $G$ son groupe adjoint. La
d\'{e}finition suivante est due \`{a} D. Panyushev.

\medskip

{\bf D\'{e}finition.}~{\em {\em (i)} On dit que deux sous-alg\`{e}bres
paraboliques $\mathfrak{p}$ et $\mathfrak{p}'$ de $\mathfrak{g}$ sont
faiblement oppos\'{e}es si $\mathfrak{p} + \mathfrak{p}' =
\mathfrak{g}$.

{\em (ii)} On appelle sous-alg\`{e}bre sp\'{e}ciale de $\mathfrak{g}$
toute sous-alg\`{e}bre de Lie $\mathfrak{q}$ de $\mathfrak{g}$ de la
forme $\mathfrak{q} = \mathfrak{p} \cap \mathfrak{p}'$, o\`{u}
$\mathfrak{p}$ et $\mathfrak{p}'$ sont des sous-alg\`{e}bres
paraboliques faiblement oppos\'{e}es de $\mathfrak{g}$.}

\medskip

{\bf 2.2.}~Soient $S,T$ des parties de $\Pi$. On v\'{e}rifie
facilement que $R_{+}^{S} \cup R_{-}^{T}$ est une partie close de
$R$. Posons :
$$
\mathfrak{p} = \mathfrak{g}^{R_{+}^{S}} \oplus \mathfrak{b}_{-} \ , \ 
\mathfrak{p}' = \mathfrak{g}^{R_{-}^{T}} \oplus \mathfrak{b}.
$$
Alors $\mathfrak{p}$ et $\mathfrak{p}'$ sont des sous-alg\`{e}bres
paraboliques faiblement oppos\'{e}es de $\mathfrak{g}$, et
$$
\mathfrak{q} = \mathfrak{p} \cap \mathfrak{p}' = \mathfrak{h} \oplus
\mathfrak{g}^{R_{+}^{S}} \oplus \mathfrak{g}^{R_{-}^{T}} 
$$
est une sous-alg\`{e}bre sp\'{e}ciale de $\mathfrak{g}$.

R\'{e}ciproquement, soient $\mathfrak{p}$ et $\mathfrak{p}'$ des
sous-alg\`{e}bres paraboliques de $\mathfrak{g}$ v\'{e}rifiant
$\mathfrak{b} \subset \mathfrak{p}$ et $\mathfrak{b}_{-} \subset
\mathfrak{p}'$. On a $\mathfrak{p} + \mathfrak{p}' = \mathfrak{g}$, et
il est imm\'{e}diat qu'il existe des parties $S$ et $T$ de $\Pi$
telles que $\mathfrak{p} \cap \mathfrak{p}' = \mathfrak{h} \oplus
\mathfrak{g}^{R_{+}^{S}} \oplus \mathfrak{g}^{R_{-}^{T}}$.

On dira qu'une sous-alg\`{e}bre sp\'{e}ciale $\mathfrak{q}$ de
$\mathfrak{g}$ est standard (relativement \`{a} $\mathfrak{h}$ et
$\Pi$) si elle s'\'{e}crit $\mathfrak{q} = \mathfrak{h} \oplus
\mathfrak{g}^{R_{+}^{S}} \oplus \mathfrak{g}^{R_{-}^{T}}$, o\`{u} $S$
et $T$ sont des parties de $\Pi$. 

\medskip

{\bf 2.3.~Proposition.}~{\em Soit $\mathfrak{q}$ une sous-alg\`{e}bre
de Lie de $\mathfrak{g}$. Les conditions suivantes sont
\'{e}quivalentes :

{\em (i)} $\mathfrak{q}$ est une sous-alg\`{e}bre sp\'{e}ciale de
$\mathfrak{g}$. 

{\em (ii)} $\mathfrak{q}$ est $G$-conjugu\'{e}e \`{a} une
sous-alg\`{e}bre sp\'{e}ciale standard de $\mathfrak{g}$.}

\smallskip

{\em Preuve.}~L'implication (ii) $\Rightarrow$ (i) est
claire. Supposons (i) v\'{e}rifi\'{e}, et prouvons (ii).

Ecrivons $\mathfrak{q} = \mathfrak{p} \cap \mathfrak{p}'$, o\`{u}
$\mathfrak{p}$ et $\mathfrak{p}'$ sont des sous-alg\`{e}bres
paraboliques faiblement oppos\'{e}es de $\mathfrak{g}$. On sait qu'il 
existe une sous-alg\`{e}bre de Cartan $\mathfrak{k}$ de $\mathfrak{g}$
contenue dans $\mathfrak{q}$. Notons $R'$ le syst\`{e}me de racines de
$(\mathfrak{g}, \mathfrak{k})$. Ecrivons $\mathfrak{p} = \mathfrak{k}
\oplus \mathfrak{g}^{P}$, $\mathfrak{p}' = \mathfrak{k} \oplus
\mathfrak{g}^{Q}$, o\`{u} $P$ et $Q$ sont des parties paraboliques de
$R'$. D'apr\`{e}s les hypoth\`{e}ses, il vient :
$$
P \cup (-P) = Q \cup (-Q) = P \cup Q = R'.
$$

On a donc aussi $(-P) \cup (-Q) = R'$. Soit $T = P \cap (-Q)$ ; c'est
une partie close de $R'$. On va montrer qu'elle est parabolique,
c'est-\`{a}-dire que $T \cup (-T) = R'$.

Soit $\alpha \in R' \backslash T$. Distinguons plusieurs cas :

a) Si $\alpha \notin P$ et $\alpha \notin -Q$, alors $\alpha \in -P$
et $\alpha \in Q$, donc $\alpha \in -T$.

b) Si $\alpha \in P$ et $\alpha \notin -Q$, on a $\alpha \in -P$ et
$\alpha \in Q$, car $(-P) \cup (-Q) = R'$.

c) Si $\alpha \notin P$ et $\alpha \in -Q$, il vient $\alpha \in -P$
et $\alpha \in Q$, car $P \cup Q = R'$.

Dans tous les cas, on a obtenu $\alpha \in -T$, ce qui fournit le
r\'{e}sultat. L'implication (i) $\Rightarrow$ (ii) est alors
imm\'{e}diate. \qed  

\bigskip

{\bf 3. Majoration de l'indice}

\bigskip

{\bf 3.1.}~Les notations des alin\'{e}as 3.1 \`{a} 3.7 seront
utilis\'{e}es dans toute la suite. On conserve les hypoth\`{e}ses de
1.3. 

Soient $\mathfrak{a}$ une sous-alg\`{e}bre de Lie de $\mathfrak{g}$ et
$X \in \mathfrak{g}$. On d\'{e}finit une forme lin\'{e}aire
$\varphi_{\mathfrak{a}}^{X}$ sur $\mathfrak{a}$ en posant, pour $Y \in
\mathfrak{a}$ :
$$
\varphi_{\mathfrak{a}}^{X}(Y) = \kappa (X,Y).
$$
Si $Z \in \mathfrak{a}$, on a $Z.\varphi_{\mathfrak{a}}^{X} =
\varphi_{\mathfrak{a}}^{[Z,X]}$ et $\mathfrak{a} \cap \mathfrak{g}^{X}
\subset \mathfrak{a}^{(\varphi_{\mathfrak{a}}^{X})}$, o\`{u}
$\mathfrak{g}^{X}$ est le centralisateur de $X$ dans $\mathfrak{g}$.

\medskip

{\bf 3.2.}~Soit $\{ H_{1}, \dots , H_{\ell}\}$ une base de
$\mathfrak{h}$, o\`{u} $\ell$ est le rang de $\mathfrak{g}$. Si
$\alpha \in R$, on note $X_{\alpha}$ un \'{e}l\'{e}ment non nul de
$\mathfrak{g}^{\alpha}$. Alors $\{ H_{i}; 1 \leqslant i \leqslant
\ell\} \cup \{ X_{\alpha}; \alpha \in R\}$ est une base de
$\mathfrak{g}$, dont on d\'{e}signe par $\{ H_{i}^{*}; 1 \leqslant i
\leqslant \ell\} \cup \{ X_{\alpha}^{*}; \alpha \in R\}$ la base
duale.

\medskip

{\bf 3.3.}~Dans la suite, on fixe un ordre total sur l'ensemble
$R_{+}$. 

Soient $S,T$ des parties de $\Pi$. On d\'{e}signe par $E_{S,T}$ le
sous-espace de $\mathfrak{h}^{*}$ en\-gen\-dr\'{e} par les
$\varepsilon_{K}$, avec $K \in \mathcal{K}(S) \cup \mathcal{K}(T)$. On
pose :
$$
\mathfrak{q} = \mathfrak{g}_{S,T} = \mathfrak{h} \oplus
\operatornamewithlimits{\textstyle \sum}_{\alpha \in R_{+}^{S}}
\mathfrak{g}^{\alpha} \oplus \operatornamewithlimits{\textstyle
\sum}_{\alpha \in R_{-}^{T}} \mathfrak{g}^{\alpha}.
$$
L'alg\`{e}bre de Lie $\mathfrak{q}$ est une sous-alg\`{e}bre
sp\'{e}ciale standard relativement \`{a} $\mathfrak{h}$ et $\Pi$.

Pour $K \in \mathcal{K}(S)$ (resp. $L  \in \mathcal{K}(T)$), on note
$a_{K}$ (resp. $b_{L}$) un \'{e}l\'{e}ment non nul de $\Bbbk$. Soit
$$
f = \operatornamewithlimits{\textstyle \sum}_{K \in \mathcal{K}(S)} a_{K}
X_{\varepsilon_{K}}^{*} + \operatornamewithlimits{\textstyle \sum}_{L
\in \mathcal{K}(T)} b_{L} X_{-\varepsilon_{L}}^{*} \in
\mathfrak{q}^{*}. 
$$
On a aussi $f = \varphi_{\mathfrak{q}}^{X_{S,T}}$, avec 
$$
X_{S,T} = \operatornamewithlimits{\textstyle
\sum}_{K \in \mathcal{K}(S)} a'_{K} 
X_{-\varepsilon_{K}} + \operatornamewithlimits{\textstyle \sum}_{L
\in \mathcal{K}(T)} b'_{L} X_{\varepsilon_{L}},
$$
o\`{u}
$$
a'_{K} = \frac{a_{K}}{\kappa (X_{\varepsilon_{K}}, X_{-\varepsilon_{K}})}
\ \raisebox{2pt}{,} \ b'_{L} = \frac{b_{L}}{\kappa (X_{\varepsilon_{L}},
X_{-\varepsilon_{L}})} \cdotp
$$

\medskip

Soient $m_{S}$ (resp. $m_{T}$) le cardinal de $\mathcal{K}(S)$
(resp. $\mathcal{K}(T)$) et $m = m_{S} + m_{T}$. On note 
$$
\Omega_{f} = (a_{K_{1}}, \dots , a_{K_{m_{S}}}, b_{L_{1}}, \dots ,
b_{m_{T}}) \in (\Bbbk \backslash \{ 0\} )^{m},
$$
avec $\varepsilon_{K_{1}} < \cdots < \varepsilon_{K_{m_{S}}}$ et
$\varepsilon_{L_{1}} < \cdots < \varepsilon_{L_{m_{T}}}$. On dit que
$\Omega_{f}$ d\'{e}finit $f$.

\medskip

{\bf 3.4.}~Pour $K \in \mathcal{K}(S)$ et $L \in \mathcal{K}(T)$, on
adopte les notations suivantes :

$\bullet$ $\mathcal{H}_{1}(K)$ est l'ensemble des couples $(\alpha ,
\varepsilon_{K} - \alpha )$, avec $\alpha \in \Gamma_{0}^{K}$ et
$\alpha < \varepsilon_{K} - \alpha$.

$\bullet$ $\mathcal{H}_{2}(L)$ est l'ensemble des couples $(-\beta ,
-\varepsilon_{L} + \beta )$, avec $\beta \in \Gamma_{0}^{L}$ et $\beta
< \varepsilon_{L} - \beta$.

$\bullet$ $\mathcal{I}_{1}(K,L)$ est l'ensemble des couples $(-\beta ,
\varepsilon_{K})$, avec $\beta \in \Gamma_{0}^{L}$, et pour lesquels
il existe $L' \in \mathcal{K}(T)$ v\'{e}rifiant $\varepsilon_{K} -
\beta = - \varepsilon_{L'}$.

$\bullet$ $\mathcal{I}_{2}(K,L)$ est l'ensemble des couples $(\alpha ,
-\varepsilon_{L})$, avec $\alpha \in \Gamma_{0}^{K}$, et pour lesquels
il existe $K' \in \mathcal{K}(T)$ v\'{e}rifiant $\alpha -
\varepsilon_{L}  =  \varepsilon_{K'}$.

$\bullet$ $\mathcal{J}_{1}(K,L)$ est l'ensemble des couples $(\alpha ,
-\beta )$, avec $\alpha \in \Gamma_{0}^{K}$ et $\beta \in
\Gamma_{0}^{L}$, et pour lesquels il existe $K' \in \mathcal{K}(S)$
v\'{e}rifiant $\alpha - \beta = \varepsilon_{K'}$.

$\bullet$ $\mathcal{J}_{2}(K,L)$ est l'ensemble des couples $(\alpha ,
-\beta )$, avec $\alpha \in \Gamma_{0}^{K}$ et $\beta \in
\Gamma_{0}^{L}$, et pour lesquels il existe $L' \in \mathcal{K}(T)$
v\'{e}rifiant $\alpha - \beta = -\varepsilon_{L'}$.

\medskip

On pose :
\begin{gather*}
\mathcal{H}_{1} = \operatornamewithlimits{\textstyle \bigcup}_{K \in
\mathcal{K}(S)} \mathcal{H}_{1}(K) \ , \ \mathcal{H}_{2} =
\operatornamewithlimits{\textstyle \bigcup}_{L \in \mathcal{K}(T)}
\mathcal{H}_{2}(L) \ , \ \mathcal{H} = \mathcal{H}_{1} \cup
\mathcal{H}_{2}, \\
\mathcal{I}_{1} = \operatornamewithlimits{\textstyle \bigcup}_{(K,L)
\in \mathcal{K}(S) \times \mathcal{K}(T)} \mathcal{I}_{1}(K,L) \ , \ 
\mathcal{I}_{2} = \operatornamewithlimits{\textstyle \bigcup}_{(K,L)
\in \mathcal{K}(S) \times \mathcal{K}(T)} \mathcal{I}_{2}(K,L), \\
\mathcal{J} = \Big(\operatornamewithlimits{\textstyle \bigcup}_{(K,L)
\in \mathcal{K}(S) \times \mathcal{K}(T)} \mathcal{J}_{1}(K,L)\Big)
\cup \Big(\operatornamewithlimits{\textstyle \bigcup}_{(K,L) \in
\mathcal{K}(S) \times \mathcal{K}(T)} \mathcal{J}_{2}(K,L)\Big), \\
\mathcal{Z} = \mathcal{H} \cup \mathcal{I}_{1} \cup \mathcal{I}_{2}
\cup \mathcal{J}, \\
r = \card \mathcal{H} \ , \ s = \dim E_{S,T}.
\end{gather*}

Si $z = (\alpha , \beta ) \in \mathcal{Z}$, on notera $\widetilde{z} =
\{ \alpha , \beta \}$ l'ensemble sous-jacent \`{a} $z$. 

\medskip

{\bf 3.5.}~{\bf Lemme.}~{\em Soient $K \in \mathcal{K}(S)$ et $L \in
\mathcal{K}(T)$. Alors :

{\em (i)} Si $(-\beta , \varepsilon_{K}) \in \mathcal{I}_{1}(K,L)$
v\'{e}rifie $\varepsilon_{K} - \beta = -\varepsilon_{L'}$, on a $L'
\subsetneq L$ et $\varepsilon_{K} \in \Gamma_{0}^{L}$.

{\em (ii)} Si $(\alpha , -\varepsilon_{L}) \in \mathcal{I}_{2}(K,L)$
v\'{e}rifie $\alpha - \varepsilon_{L} = \varepsilon_{K'}$, on a $K'
\subsetneq K$ et $\varepsilon_{L} \in \Gamma_{0}^{K}$.

{\em (iii)} Si $(\alpha , -\beta ) \in \mathcal{J}_{1}(K,L)$
v\'{e}rifie $\alpha - \beta = \varepsilon_{K'}$, on a $K' \subsetneq
K$ et $\beta \in \Gamma_{0}^{K}$.

{\em (iv)} Si $(\alpha , -\beta ) \in \mathcal{J}_{2}(K,L)$ v\'{e}rifie
$\alpha - \beta = -\varepsilon_{L'}$, on a $L' \subsetneq L$ et
$\alpha \in \Gamma_{0}^{L}$.}

\smallskip

{\em Preuve.}~C'est imm\'{e}diat d'apr\`{e}s les lemmes 1.7 et
1.8. \qed 

\medskip

{\bf 3.6.}~Si $K \in \mathcal{K}(S)$ et $L \in \mathcal{K}(T)$, on
note :
\begin{gather*}
\mathfrak{a}_{K}^{\bullet} = \operatornamewithlimits{\textstyle
\sum}_{\alpha \in \Gamma_{0}^{K}} \mathfrak{g}^{\alpha} \ , \
\mathfrak{b}_{L}^{\bullet} = \operatornamewithlimits{\textstyle
\sum}_{\alpha \in \Gamma_{0}^{L}} \mathfrak{g}^{-\alpha}, \\
\mathfrak{a}_{K} = \Bbbk X_{\varepsilon_{K}} +
\mathfrak{a}_{K}^{\bullet} \ , \ \mathfrak{b}_{L} = \Bbbk
X_{-\varepsilon_{L}} + \mathfrak{b}_{L}^{\bullet}, \\
\mathfrak{u} = \operatornamewithlimits{\textstyle \sum}_{K \in
\mathcal{K}(S)} \Bbbk X_{\varepsilon_{K}} +
\operatornamewithlimits{\textstyle \sum}_{L \in \mathcal{K}(T)} \Bbbk
X_{-\varepsilon_{L}} \ , \ \mathfrak{v} =
\operatornamewithlimits{\textstyle \sum}_{K \in \mathcal{K}(S)}
\mathfrak{a}_{K}^{\bullet} + \operatornamewithlimits{\textstyle
\sum}_{L \in \mathcal{K}(T)} \mathfrak{b}_{L}^{\bullet}.
\end{gather*}

D'apr\`{e}s 1.8, $\mathfrak{a}_{K}$ et $\mathfrak{b}_{L}$ sont des
alg\`{e}bres de Heisenberg, de centres respectifs $\Bbbk
X_{\varepsilon_{K}}$ et $\Bbbk X_{-\varepsilon_{L}}$.

La dimension de $\mathfrak{a}_{K}^{\bullet}$
(resp. $\mathfrak{b}_{L}^{\bullet}$) est un entier pair $2n_{K}$
(resp. $2n_{L}$), et :
$$
r = \operatornamewithlimits{\textstyle \sum}_{K \in \mathcal{K}(S)}
n_{K} + \operatornamewithlimits{\textstyle \sum}_{L \in
\mathcal{K}(T)} n_{L}.
$$

On identifie $\mathfrak{h}^{*}$, $\mathfrak{u}^{*}$ et
$\mathfrak{v}^{*}$ \`{a} des sous-espaces de $\mathfrak{q}^{*}$ au
moyen de la d\'{e}\-com\-po\-si\-tion $\mathfrak{q} = \mathfrak{h} \oplus
\mathfrak{u} \oplus \mathfrak{v}$.

\medskip

{\bf 3.7.}~En utilisant la notation de 3.2, pour $z = (\alpha , \beta )
\in \mathcal{Z}$, on note :
$$
v_{z} = X_{\alpha}^{*} \wedge X_{\beta}^{*} \in {\textstyle
\bigwedge}^{2} \mathfrak{q}^{*}.
$$

Soient $K,K' \in \mathcal{K}(S)$, $L,L' \in
\mathcal{K}(T)$, $\alpha \in \Gamma_{0}^{K}$, et $\beta \in
\Gamma_{0}^{L}$. Compte tenu des lemmes 1.7 et 1.8, on a les
r\'{e}sultats suivants :

$\bullet$ $\varepsilon_{K} - \varepsilon_{L} \notin \{
\varepsilon_{K'}, - \varepsilon_{L'}\}$.

$\bullet$ $\alpha - \varepsilon_{L} \ne -\varepsilon_{L'}$.

$\bullet$ $\varepsilon_{K} - \beta \ne \varepsilon_{K'}$.

Il en r\'{e}sulte que, si l'on identifie la forme bilin\'{e}aire
altern\'{e}e $\Phi_{f}$ sur $\mathfrak{q}$ \`{a} un \'{e}l\'{e}ment de
$\bigwedge^{2} \mathfrak{q}^{*}$, alors 
$$
\Phi_{f} = \Psi_{f} + \Theta_{f},
$$
avec $\Theta_{f} \in E_{S,T} \bigwedge \mathfrak{u}^{*}$ et
$$
\Psi_{f} = \operatornamewithlimits{\textstyle \sum}_{z \in
\mathcal{Z}} \lambda_{z} v_{z},
$$
avec $\lambda_{z} \in \Bbbk$ pour tout $z \in
\mathcal{Z}$.

Si $z = (\alpha , \beta ) \in \mathcal{Z}$, on a
$\Theta_{f}(X_{\alpha}, X_{\beta}) = 0$. D'autre part, ou
$[X_{\alpha}, X_{\beta}] = \mu_{z}X_{\varepsilon_{K}}$, avec $K \in
\mathcal{K}(S)$, ou $[X_{\alpha}, X_{\beta}] =
\mu_{z}X_{-\varepsilon_{L}}$, avec $L \in \mathcal{K}(T)$, le scalaire
$\mu_{z}$ \'{e}tant non nul. Par cons\'{e}quent :
$$
\lambda_{z} = \Phi_{f} (X_{\alpha},
X_{\beta}) = f([X_{\alpha}, X_{\beta}]) = \begin{cases}
\mu_{z} a_{K} \text{ si } \ [X_{\alpha}, X_{\beta}] = \mu_{z}
X_{\varepsilon_{K}} ,\hfill\\
\mu_{z} b_{L} \text{ si } \ [X_{\alpha}, X_{\beta}] = \mu_{z}
X_{-\varepsilon_{L}}. \hfill
\end{cases}
$$
On voit donc que $\lambda_{z}$ est non nul.

\medskip

{\bf 3.8.}~{\bf Lemme.}~{\em On conserve les notations
pr\'{e}c\'{e}dentes. Alors :

{\em (i)}~$\mathfrak{q}^{(f)}$ contient une sous-alg\`{e}bre
commutative de $\mathfrak{g}$, form\'{e}e d'\'{e}l\'{e}ments
semi-simples, et de dimension :
$$
\dim \mathfrak{h} - \dim E_{S,T} + \card \big( \mathcal{K}(S) \cap
\mathcal{K}(T)\big) .
$$

{\em (ii)}~On a $\wedge^{s} \Theta_{f} \ne 0$ et $\wedge^{s+1}
\Theta_{f} = 0$. 

{\em (iii)}~Il existe un ouvert non vide $U$ de $(\Bbbk \backslash \{
0\} )^{m}$ tel que $\wedge^{r} \Psi_{f} \ne 0$ d\`{e}s que
$\Omega_{f} \in U$.

{\em (iv)}~Si $\Omega_{f} \in U$, alors $\wedge^{r+s} \Phi_{f} \ne
0$.}

\smallskip

{\em Preuve.}~(i) L'orthogonal $\mathfrak{k}$ de $E_{S,T}$ dans
$\mathfrak{h}$ est de dimension $\dim \mathfrak{h} - \dim
E_{S,T}$. Avec la notation $X_{S,T}$ de 3.3, si $H \in \mathfrak{k}$,
on a $[H,X_{S,T}] = 0$, soit $H.f = 0$.

Pour $K \in \mathcal{K}(S)\cap \mathcal{K}(T)$, posons :
$$
Y_{K} = a'_{K} X_{\varepsilon_{K}} + b'_{K} X_{-\varepsilon_{K}}.
$$
Il est bien connu que $Y_{K}$ est un \'{e}l\'{e}ment semi-simple de
$\mathfrak{g}$. D'autre part, d'apr\`{e}s 3.3 et l'assertion (iii) du
lemme 1.7, on a $Y_{K}.f = 0$.

Compte tenu de ces remarques, l'assertion est imm\'{e}diate.

(ii)~Si $k, l$ sont des entiers strictement positifs, nous noterons
$\mathfrak{M}_{k,l}$ l'ensemble des matrices \`{a} $k$ lignes et
$l$ colonnes \`{a} \'{e}l\'{e}ments dans $\Bbbk$. Pour $A \in
\mathfrak{M}_{k,l}$, ${}^{t}A$ d\'{e}signe la transpos\'{e}e de
$A$. 

Soit $(\alpha_{1}, \dots , \alpha_{s})$ une base de $E_{S,T}$ telle
que $\alpha_{i} = \varepsilon_{K_{i}}$ si $1 \leqslant i \leqslant p$,
et $\alpha_{i} = -\varepsilon_{L_{i}}$ si $p+1 \leqslant i \leqslant
s$, les $K_{i}$ \'{e}tant des \'{e}l\'{e}ments de $\mathcal{K}(S)$ et
les $L_{i}$ des \'{e}l\'{e}ments de $\mathcal{K}(T)$. Compl\'{e}tons
cette base en une base $\mathcal{B}' = (\alpha_{1}, \dots ,
\alpha_{\ell})$ de $\mathfrak{h}^{*}$, et soit $\mathcal{B} = (h_{1},
\dots , h_{\ell})$ la base de $\mathfrak{h}$ duale de
$\mathcal{B}'$. De m\^{e}me, compl\'{e}tons le syst\`{e}me libre
$(X_{\varepsilon_{K_{1}}}, \dots , X_{\varepsilon_{K_{p}}},
X_{-\varepsilon_{L_{p+1}}}, \dots , X_{-\varepsilon_{L_{s}}})$ en une
base $\mathcal{C}$ de $\mathfrak{u}$. Alors $\mathcal{D} = \mathcal{B}
\cup \mathcal{C}$ est une base de $\mathfrak{h} + \mathfrak{u}$. Si
$s+1 \leqslant k \leqslant \ell$, on a $[h_{k}, \mathfrak{u}] = \{ 0\}$.

Avec les notations usuelles, les racines $\varepsilon_{K}$
(resp. $\varepsilon_{L}$) sont fortement orthogonales, et
$\varepsilon_{K} - \varepsilon_{L} \notin \mathcal{K}(S) \cup
\mathcal{K}(T)$. Il en r\'{e}sulte que, dans la base $\mathcal{D}$, la
restriction de la forme bilin\'{e}aire $\Phi_{f}$ (qui est aussi la
restriction de $\Theta_{f}$) a une matrice de la forme 
$$
M = \begin{pmatrix}
0 & 0 & D & A\\
0 & 0 & 0 & 0\\
-{}^{t}D & 0 & 0 & 0\\
-{}^{t}A & 0 & 0 & 0\end{pmatrix},
$$
o\`{u} $A \in \mathfrak{M}_{s,m-s}$ ($m$ est d\'{e}fini en 3.3), et
o\`{u} $D \in \mathfrak{M}_{s,s}$ est diagonale, sa diagonale
\'{e}tant :
$$
(a_{K_{1}}, \dots , a_{K_{p}}, b_{L_{p+1}}, \dots , b_{L_{s}}).
$$
On en d\'{e}duit que $M$ est de rang $2s$. Il est alors bien connu que
que $\wedge^{s} \Theta_{f} \ne 0$ et que $\wedge^{s+1} \Theta_{f} =
0$.

(iii) Notons $z_{1}, \dots , z_{n}$ les \'{e}l\'{e}ments de
$\mathcal{Z}$, en convenant que $z_{1}, \dots , z_{r}$ sont ceux de
$\mathcal{H}$. Afin de simplifier les notations, on \'{e}crira
$\lambda_{i} v_{i}$ pour $\lambda_{z_{i}} v_{z_{i}}$ et $\mu_{i}$ pour
$\mu_{z_{i}}$.

Si $z,z' \in \mathcal{Z}$, on a $v_{z} \wedge v_{z'} = v_{z'} \wedge v_{z}$ et $v_{z}
\wedge v_{z} = 0$. Par cons\'{e}quent :
$$
\wedge^{r} \Psi_{f} = r! \operatornamewithlimits{\textstyle \sum}_{1
\leqslant i_{1} < \cdots < i_{r} \leqslant n} \lambda_{i_{1}} \cdots
\lambda_{i_{r}} v_{i_{1}} \wedge \cdots \wedge v_{i_{r}}. 
$$

Dans la somme pr\'{e}c\'{e}dente, le coefficient de $v_{1} \wedge
\cdots \wedge v_{r}$ est :
$$
\operatornamewithlimits{\textstyle \prod}_{K \in \mathcal{K}(S)}
a_{K}^{n_{K}} \Big(\operatornamewithlimits{\textstyle \prod}_{z \in
\mathcal{H}_{1}(K)} \mu_{z}\Big) \operatornamewithlimits{\textstyle
\prod}_{L \in \mathcal{K}(T)} b_{L}^{n_{L}}
\Big(\operatornamewithlimits{\textstyle \prod}_{z \in
\mathcal{H}_{2}(L)} \mu_{z}\Big).
$$

Supposons que $v_{i_{1}} \wedge \cdots \wedge v_{i_{r}} = \lambda
v_{1} \wedge \cdots \wedge v_{r}$, avec $\lambda \in \Bbbk \backslash
\{ 0\}$, o\`{u} $i_{1} < \cdots < i_{r}$ et $(i_{1}, \dots , i_{r})
\ne (1, \dots , r)$. L'ensemble $\mathcal{S} = \widetilde{z}_{1} \cup
\cdots \cup \widetilde{z_{r}}$ est alors la r\'{e}union disjointe des
ensembles $\widetilde{z_{i_{t}}}$ pour $1 \leqslant t \leqslant r$. On
en d\'{e}duit que, si $z_{i_{t}} \notin \mathcal{H}$, on a $z_{i_{t}}
\in \mathcal{J}$, et il existe donc $K \in \mathcal{K}(S)$ et $L \in
\mathcal{K}(T)$ tels que $\Gamma_{0}^{K} \cap \widetilde{z_{i_{t}}}
\ne \emptyset$ et $(-\Gamma_{0}^{L})\cap \widetilde{z_{i_{t}}} \ne
\emptyset$. Notons $z_{j_{1}}, \dots , z_{j_{k}}$ les \'{e}l\'{e}ments
$z_{i_{t}}$ n'appartenant pas \`{a} $\mathcal{H}$.

Fixons $K_{0} \in \mathcal{K}(S)$ maximal pour l'inclusion parmi les
\'{e}l\'{e}ments $K$ de $\mathcal{K}(S)$ v\'{e}rifiant :
$$
\Gamma_{0}^{K} \cap (\widetilde{z_{j_{1}}} \cup \cdots \cup
\widetilde{z_{j_{k}}}) \ne \emptyset .
$$ 

Il existe $t \in \{ 1, \dots , k\}$ tel que, si $z = z_{_{j_{t}}}$, on
ait $\Gamma_{0}^{K_{0}} \cap \widetilde{z} \ne \emptyset$. Soit
$\alpha \in \Gamma_{0}^{K_{0}} \cap \widetilde{z}$.

Il en r\'{e}sulte que $z \in \mathcal{J}$ et on a $(\alpha ,
\varepsilon_{K_{0}} - \alpha ) \ne z_{i_{l}}$ pour $1 \leqslant l
\leqslant r$. D'apr\`{e}s la fin de 3.7 et le lemme 3.5, il vient ou
$\lambda_{z} = \mu_{z}a_{K}$, avec $K \in \mathcal{K}(S)$ et $K \ne
K_{0}$, ou $\lambda_{z} = \mu_{z} b_{L}$, avec $L \in \mathcal{K}(T)$.

On d\'{e}duit de ceci que le coefficient de $v_{i_{1}} \wedge \cdots
\wedge v_{i_{r}}$ dans la somme donnant $\wedge^{r} \Psi_{f}$ est
de la forme 
$$
\mu_{i_{1}} \cdots \mu_{i_{r}} \operatornamewithlimits{\textstyle
\prod}_{K \in \mathcal{K}(S)} a_{K}^{m_{K}}
\operatornamewithlimits{\textstyle \prod}_{L \in \mathcal{K}(T)}
b_{L}^{m_{L}},
$$
avec $m_{K_{0}} < n_{K_{0}}$. 

Il est alors clair qu'il existe un ouvert non vide $U$ de $(\Bbbk
\backslash \{ 0\} )^{m}$ v\'{e}rifiant $\wedge^{r} \Psi_{f} \ne 0$
si $\Omega_{f} \in U$. 

(iv)~On a :
$$
\wedge^{r+s} \Phi_{f} = \operatornamewithlimits{\textstyle
\sum}_{k=0}^{r+s} \begin{pmatrix}
r+s \\
k\end{pmatrix} (\wedge^{k} \Psi_{f} ) \wedge (\wedge^{r+s-k}
\Theta_{f}). 
$$

Comme $\wedge^{j}\Theta_{f} \in (\textstyle{\bigwedge^{j}}
E_{S,T})\bigwedge (\textstyle{\bigwedge^{j}}\mathfrak{u}^{*})$, pour
montrer que $\wedge^{r+s} \Phi_{f} \ne 0$, il suffit de prouver que
$(\wedge^{r} \Psi_{f})\wedge (\wedge^{s} \Theta_{f}) \ne 0$.

Conservons les notations $v_{1}, \dots , v_{r}$ et $\mathcal{S}$ de
(iii). On a vu que, si $\Omega_{f} \in U$,
$$
\wedge^{r} \Psi_{f} = \lambda v_{1}\wedge \cdots \wedge v_{r} + w,
$$
o\`{u} $\lambda \in \Bbbk \backslash \{ 0\}$ et o\`{u} $w$ est une
combinaison lin\'{e}aire de termes $v_{z_{i_{1}}} \wedge \cdots \wedge
v_{z_{i_{r}}}$, avec $\widetilde{z_{i_{1}}} \cup \cdots \cup
\widetilde{z_{i_{r}}} \ne \mathcal{S}$. Il est donc imm\'{e}diat que
$(\wedge^{r} \Psi_{f}) \wedge (\wedge^{s} \Theta_{f}) \ne 0$ si
$\Omega_{f} \in U$. \qed

\medskip 

{\bf Remarque.}~On voit facilement que la preuve pr\'{e}c\'{e}dente
montre que, pour $\Omega_{f} \in U$, la restriction de $\Phi_{f}$
\`{a} $\mathfrak{v} \times \mathfrak{v}$ est non
d\'{e}g\'{e}n\'{e}r\'{e}e. 

\medskip

{\bf 3.9.}~{\bf Th\'{e}or\`{e}me.}~{\em Avec les notations
pr\'{e}c\'{e}dentes, on a :
$$
\chi (\mathfrak{g}_{S,T}) \leqslant \rg (\mathfrak{g}) + \card
\mathcal{K}(S) + \card \mathcal{K}(T) - 2 \dim E_{S,T}.
$$}

\smallskip

{\em Preuve.}~De $\mathfrak{q} = \mathfrak{h} \oplus \mathfrak{u}
\oplus \mathfrak{v}$, on d\'{e}duit que :
$$
\dim \mathfrak{g} = \dim \mathfrak{h} + \card \mathcal{K}(S) + \card
\mathcal{K}(T) + 2r.
$$

Si $\Omega_{f} \in U$, le fait que $\wedge^{r+s} \Phi_{f} \ne 0$
signifie que $\rg (\Phi_{f}) \geqslant 2(r+s)$, c'est-\`{a}-dire :
$$
\dim \mathfrak{q}^{(f)} \leqslant \dim \mathfrak{g} - 2 (r+s).
$$
Il vient donc :
$$
\dim \mathfrak{q}^{(f)} \leqslant \dim \mathfrak{h} + \card
\mathcal{K}(S) + \card \mathcal{K}(T) - 2s.
$$
D'o\`{u} le r\'{e}sultat. \qed

\medskip

{\bf 3.10.}~{\bf Remarque.}~On a $\card \mathcal{K}(S) = \dim E_{S}$,
$\card \mathcal{K}(T) = E_{T}$, ainsi que $E_{S,T} = E_{S} +
E_{T}$. L'in\'{e}galit\'{e} pr\'{e}c\'{e}dente s'\'{e}crit donc encore
:
$$
\chi (\mathfrak{g}_{S,T}) \leqslant \rg (\mathfrak{g}) + \dim E_{S} +
\dim E_{T} - 2 \dim (E_{S} + E_{T}).
$$

\medskip

{\bf 3.11.}~{\bf Corollaire.}~{\em Soit $\mathfrak{q}$ une
sous-alg\`{e}bre sp\'{e}ciale de $\mathfrak{g}$. Alors :

{\em (i)} On a $\chi (\mathfrak{q}) \leqslant \rg (\mathfrak{g})$.

{\em (ii)} Dire que $\chi (\mathfrak{q}) = \rg (\mathfrak{g})$
signifie que $\mathfrak{q}$ est une sous-alg\`{e}bre de Levi de
$\mathfrak{g}$. }

\smallskip

{\em Preuve.}~La proposition 2.3 montre que l'on peut supposer
$\mathfrak{q} = \mathfrak{g}_{S,T}$. 

(i)~D'apr\`{e}s 3.10, on a $\chi(\mathfrak{q}) \leqslant \rg
(\mathfrak{g})$.

(ii)~Si $\mathfrak{q}$ est une sous-alg\`{e}bre de Levi de
$\mathfrak{g}$, il est bien connu que $\chi (\mathfrak{q}) = \rg
(\mathfrak{g})$. 

Supposons $\chi(\mathfrak{q}) = \rg (\mathfrak{g})$. D'apr\`{e}s 3.10,
il vient $E_{S} = E_{T}$. Soit $L$ une composante connexe de $T$. On a
$\varepsilon_{L} \in E_{T} = E_{S} \subset \mathbb{Z}S$. Comme
$\varepsilon_{L}$ est une combinaison lin\'{e}aire \`{a} coefficients
entiers strictement positifs de toutes les racines de $L$ et que $\Pi$
est une base de $\mathfrak{h}^{*}$, on obtient $L \subset S$. Par
suite $T \subset S$, et $S = T$ en \'{e}changeant les r\^{o}les de $S$
et de $T$. Ainsi, $\mathfrak{q}$ est une sous-alg\`{e}bre de Levi de
$\mathfrak{g}$. \qed

\medskip

{\bf Remarque.}~Le corollaire 3.11 est une r\'{e}ponse positive \`{a}
une conjecture de D. Panyushev formul\'{e}e dans [7].

\medskip

{\bf 3.12.}~{\bf Corollaire.}~{\em Si l'ensemble des
$\varepsilon_{M}$, avec $M \in \mathcal{K}(S) \cup \mathcal{K}(T)$, est
une base de $\mathfrak{h}^{*}$, alors $\chi (\mathfrak{g}_{S,T}) =
0$.}

\smallskip

{\em Preuve.}~Si les hypoth\`{e}ses du corollaire sont
v\'{e}rifi\'{e}es, l'in\'{e}galit\'{e} de 3.10 montre que $\chi
(\mathfrak{g}_{S,T}) \leqslant 0$. \qed

\bigskip

{\bf 4.~Etude de cas particuliers et une conjecture}

\bigskip

{\bf 4.1.}~Avec les notations du paragraphe 3, on pose :
\begin{align*}
d_{S,T} & = \rg (\mathfrak{g}) + \card \mathcal{K}(S) + \card
\mathcal{K}(T) -2 \dim E_{S,T} \\
{} & = \rg (\mathfrak{g}) + \dim E_{S} + \dim E_{T} - 2 \dim (E_{S} +
E_{T}). 
\end{align*}
D'apr\`{e}s 3.9, on a $\chi (\mathfrak{g}_{S,T}) \leqslant
d_{S,T}$. On va montrer que, dans certains cas, on peut affirmer que
$\chi (\mathfrak{g}_{S,T}) = d_{S,T}$.

\medskip

{\bf 4.2}~{\bf Proposition.}~{\em Si $d_{S,T} \in \{ 0,1\}$, alors
$\chi (\mathfrak{g}_{S,T}) = d_{S,T}$.}

\smallskip

{\em Preuve.}~Le cas $d_{S,T} = 0$ est clair d'apr\`{e}s
3.9. Supposons $d_{S,T} = 1$.

On a $\dim \mathfrak{g}_{S,T} - d_{S,T} = 2(r+s)$. Comme $\dim
\mathfrak{g}_{S,T} - \chi (\mathfrak{g}_{S,T})$ est un entier pair
(c'est le rang d'une forme bilin\'{e}aire altern\'{e}e sur
$\mathfrak{g}_{S,T}$), on voit que $d_{S,T}$ et $\chi
(\mathfrak{g}_{S,T})$ ont m\^{e}me parit\'{e}. D'o\`{u} $\chi
(\mathfrak{g}_{S,T}) = 1$. \qed

\medskip

{\bf 4.3.}~Soit $\mathfrak{a}$ une sous-alg\`{e}bre de Lie
alg\'{e}brique de $\mathfrak{g}$, $A$ son groupe adjoint, et $g \in
\mathfrak{a}^{*}$. On dit que $g$ est stable s'il existe un ouvert $U$
de $\mathfrak{a}^{*}$ contenant $g$ tel que $\mathfrak{a}^{(g)}$ et
$\mathfrak{a}^{(h)}$ soient $A$-conjugu\'{e}s pour tout $h \in U$. 

Le r\'{e}sultat suivant est prouv\'{e} dans [9] :

\smallskip

{\bf Proposition.}~{\em Soient $\mathfrak{a}$ une sous-alg\`{e}bre de
Lie alg\'{e}brique de $\mathfrak{g}$ et $g \in
\mathfrak{a}^{*}$. Alors :

{\em (i)} Si $g$ est stable, c'est un \'{e}l\'{e}ment r\'{e}gulier de
$\mathfrak{a}^{*}$.

{\em (ii)} Si $\mathfrak{g}^{(f)}$ est commutative et
compos\'{e}e d'\'{e}l\'{e}ments semi-simples de $\mathfrak{g}$, c'est
un \'{e}l\'{e}ment stable de $\mathfrak{a}^{*}$.}

\medskip

{\bf 4.4.}~{\bf Proposition.}~{\em On suppose que les parties $S$ et
$T$ de $\Pi$ v\'{e}rifient l'une ou l'autre des propri\'{e}t\'{e}s
suivantes :

{\em (i)} $\mathcal{K}(S) \subset \mathcal{K}(T)$.

{\em (ii)} $\mathcal{K}(T) \subset \mathcal{K}(S)$.

{\em (iii)} La famille $(\varepsilon_{K})_{K \in \mathcal{K}(S)} \cup
(\varepsilon_{L})_{L \in \mathcal{K}(T)}$ est libre (c'est le cas en
particulier si $S \cap T = \emptyset$).

Alors, si $\Omega_{f} \in U$, $f$ est un \'{e}l\'{e}ment stable de
$\mathfrak{g}_{S,T}^{*}$, et $\chi (\mathfrak{g}_{S,T}) = d_{S,T}$.}

\smallskip

{\em Preuve.}~Si (i) est r\'{e}alis\'{e}, alors $\mathcal{K}(S) =
\mathcal{K}(S) \cap \mathcal{K}(T)$ et $E_{S,T} = E_{T}$. De m\^{e}me,
la condition (ii) implique $\mathcal{K}(T) = \mathcal{K}(S) \cap
\mathcal{K}(T)$ et $E_{S,T} = E_{S}$. Enfin, si l'hypoth\`{e}se (iii)
est vraie, on a $E_{S,T} = E_{S} \oplus E_{T}$ et $\mathcal{K}(S) \cap
\mathcal{K}(T) = \emptyset$. 

Dans tous les cas, on obtient :
$$
d_{S,T} = \rg (\mathfrak{g}) - \dim E_{S,T} + \card \big(
\mathcal{K}(S) \cap \mathcal{K}(T)\big) .
$$

Compte tenu de 3.8 et 3.9, on voit que, si $\Omega_{f} \in U$, alors $\dim
\mathfrak{q}^{(f)} = d_{S,T}$, et que l'alg\`{e}bre de Lie
$\mathfrak{q}^{(f)}$ est form\'{e}e d'\'{e}l\'{e}ments semi-simples de
$\mathfrak{g}$. On a donc le r\'{e}sultat d'apr\`{e}s 4.3. \qed

\medskip

{\bf Remarque.}~En prenant $S = \Pi$ et $T = \emptyset$, la
proposition pr\'{e}c\'{e}dente permet de retrouver le r\'{e}sultat
suivant de [9] : il existe une forme lin\'{e}aire stable dans
$\mathfrak{b}^{*}$, et on a :
$$
\chi (\mathfrak{b}) = \rg (\mathfrak{g}) - \card \mathcal{K}(\Pi ).
$$

\medskip

{\bf 4.5.}~Dans la suite, on va utiliser le r\'{e}sultat suivant ([3],
lemma 1.12.2) :

\medskip

{\bf Lemme.}~{\em Soient $V$ un espace vectoriel de dimension finie,
$V'$ un hyperplan de $V$, $\Phi$ une forme bilin\'{e}aire altern\'{e}e
sur $V$, et $\Phi'$ la restriction de $\Phi$ \`{a} $V'$. On note $N$
et $N'$ les noyaux de $\Phi$ et $\Phi'$. 

{\em (i)} Si $N \subset N'$, alors $N$ est un hyperplan de $N'$.

{\em (ii)} Si $N \not\subset N'$, on a $N' = N \cap V'$, et $N'$ est
un hyperplan de $N$.}

\medskip

{\bf 4.6.}~{\bf Proposition.}~{\em On suppose que $\card (S) = 1$ ou
que $\card (T) = 1$. Alors $\chi (\mathfrak{g}_{S,T}) = d_{S,T}$. }

\smallskip

{\em Preuve.}~Traitons le cas o\`{u} $\card (T) = 1$, soit $T = \{
\alpha \}$. On note toujours $\mathfrak{q}$ pour
$\mathfrak{g}_{S,T}$. 

Compte tenu de 4.4, on peut supposer $\alpha \in S$, $\{ \alpha\}
\notin \mathcal{K}(S)$, et $\alpha \in E_{S}$. Posons :
$$
\mathfrak{a} = \mathfrak{g}_{S, \emptyset} = \mathfrak{h} \oplus
\mathfrak{g}^{R_{+}^{S}}. 
$$

Soient
$$
U = \operatornamewithlimits{\textstyle \sum}_{K \in \mathcal{K}(S)}
a_{K} X_{-\varepsilon_{K}} \ , V = X_{\alpha} + U \ , \
f = \varphi_{\mathfrak{a}}^{U} \ , \ g = \varphi_{\mathfrak{q}}^{V},
$$
o\`{u} les coefficients $a_{K}$ sont choisis pour que $f$ soit un
\'{e}l\'{e}ment stable de $\mathfrak{a}^{*}$ (c'est possible
d'apr\`{e}s 4.4). La restriction de $\Phi_{g}$ \`{a} $\mathfrak{a}$
est $\Phi_{f}$. La preuve de 3.8, (i) montre que :
$$
\mathfrak{a}^{(f)} = \operatornamewithlimits{\textstyle \bigcap}_{K
\in \mathcal{K}(S)} \ker \varepsilon_{K}.
$$

De $\alpha \in E_{S}$, on d\'{e}duit que $[X_{\alpha},
\mathfrak{a}^{(f)}] = \{ 0\}$, soit $\mathfrak{a}^{(f)} \subset
\mathfrak{q}^{(g)}$. D'apr\`{e}s le lemme 4.5, il vient 
alors 
$$
\mathfrak{q}^{(g)} = \Bbbk  X \oplus \mathfrak{a}^{(f)},
$$
avec $X = X_{-\alpha} + Y$, o\`{u} $Y \in \mathfrak{a}$.

Comme $\mathfrak{a}^{(f)} \subset \mathfrak{h}$, on a $[\mathfrak{a},
\mathfrak{a}^{(f)}] \cap \mathfrak{a}^{(f)} = \{ 0\}$, puis $[X,
\mathfrak{a}^{(f)}] = \{ 0\}$. Par cons\'{e}quent, l'alg\`{e}bre de
Lie $\mathfrak{q}^{(g)}$ est commutative.

Soient $G$ le groupe adjoint de $\mathfrak{g}$ et $A$ le plus petit
sous-groupe alg\'{e}brique de $G$ d'alg\`{e}bre de Lie
$\ad_{\mathfrak{g}} \mathfrak{a}$. L'ensemble des restrictions des
\'{e}l\'{e}ments de $A$ \`{a} $\mathfrak{a}$ s'identifie au groupe
adjoint de $\mathfrak{a}$. 

L'ensemble des \'{e}l\'{e}ments stables de $\mathfrak{a}^{*}$ est un
ouvert non vide $A$-invariant. On en d\'{e}duit qu'il existe $h \in
\mathfrak{q}^{*}$ r\'{e}gulier tel que $\lambda = h | \mathfrak{a}$
soit stable. Si $\theta \in A$, on a $\theta (h)|\mathfrak{a} = \theta
(\lambda )$. Par cons\'{e}quent, on peut supposer que :
$$
\mathfrak{a}^{(\lambda )} = \mathfrak{a}^{(f)} =
\operatornamewithlimits{\textstyle \bigcap}_{K \in \mathcal{K}(S)} \ker
\varepsilon_{K} \subset \mathfrak{h}.
$$

On a $h ([\mathfrak{a}, \mathfrak{a}^{(\lambda )}]) = \{ 0\}$. D'autre
part, $X_{-\alpha}$ commute \`{a} $\mathfrak{a}^{(f)} =
\mathfrak{a}^{(\lambda )}$. On a donc $h([X_{-\alpha},
\mathfrak{a}^{(\lambda )}]) = \{ 0\}$. Par suite, $h([\mathfrak{q},
\mathfrak{a}^{(\lambda )}]) = \{ 0\}$ et $\mathfrak{a}^{(\lambda )}
\subset \mathfrak{q}^{(h)}$. 

D'apr\`{e}s 4.5, on obtient $\dim \mathfrak{q}^{(h)} = 1 + \dim
\mathfrak{a}^{(\lambda )}$, d'o\`{u} $\chi (\mathfrak{q}) = 1 + \chi
(\mathfrak{a})$.

On a ici $E_{S,T} = E_{S,\emptyset}$, $\card \mathcal{K}(T) = 1$, et
$\chi (\mathfrak{a}) = d_{S,\emptyset}$. On en d\'{e}duit que $\chi
(\mathfrak{q}) = d_{S,T}$. \qed

\medskip

{\bf Remarque.}~Il est donn\'{e} dans [9] un exemple d'alg\`{e}bre de
Lie v\'{e}rifiant les conditions de 4.6, mais ne poss\'{e}dant aucune
forme lin\'{e}aire stable. 

\medskip

{\bf 4.7.}~Compte tenu des exemples trait\'{e}s pr\'{e}c\'{e}demment,
nous \'{e}noncerons l'assertion suivante :

\medskip

{\bf Conjecture.}~{\em Si $S$ et $T$ sont des parties de $\Pi$, on a
$\chi (\mathfrak{g}_{S,T}) = d_{S,T}$.}

\medskip

{\bf 4.8.}~Faisons quelques remarques.

1) Supposons $\mathfrak{g}$ de type $A$. L'indice d'une
sous-alg\`{e}bre sp\'{e}ciale de $\mathfrak{g}$ est calcul\'{e} dans
[2] au moyen d'une formule combinatoire. On peut v\'{e}rifier qu'elle
confirme la conjecture 4.7.

2) Pour $\mathfrak{g}$ de type $A,B,C$ ou $D$, des formules de
r\'{e}currence reliant l'indice des sous-alg\`{e}bres sp\'{e}ciales de
$\mathfrak{g}$ sont \'{e}tablies dans [7]. Ces formules sont compatibles
avec 4.7.

3) Des calculs effectu\'{e}s par D. Panyushev et R. Ushirobira lorsque
$\mathfrak{g}$ est de type $F_{4}$ ou $G_{2}$ confirment la
conjecture 4.7.

4) D'apr\`{e}s 4.2, 4.4, et 4.6, la conjecture 4.7 est
 v\'{e}rifi\'{e}e si $\rg (\mathfrak{g}) \leqslant 2$.

\bigskip

{\bf 5.~Sur l'indice des sous-alg\`{e}bres paraboliques}

\bigskip

{\bf 5.1.}~{\bf Lemme.}~{\em Soient $S$ une partie de $\Pi$ et $n =
\card \mathcal{K}(S)$. Il existe des parties $S_{1}, \dots , S_{n}$ de
$\Pi$ v\'{e}rifiant les conditions suivantes :

{\em (i)} $S_{1} \subset S_{2} \subset \cdots \subset S_{n} = S$.

{\em (ii)} $\mathcal{K}(S_{1}) \subset \mathcal{K}(S_{2}) \subset
\cdots \subset \mathcal{K}(S_{n})$.

{\em (iii)} $\card \mathcal{K}(S_{i}) = i$ pour $1 \leqslant i
\leqslant n$.}

\smallskip

{\em Preuve.}~Soient $K$ une composante connexe de $S$ et $K'$
l'ensemble des $\alpha \in K$ qui v\'{e}rifient $\langle \alpha ,
\varepsilon_{K}^{\vee}\rangle = 0$. Si $S_{n-1} = S \backslash (K
\backslash K')$, il est clair, d'apr\`{e}s 1.5, que $\card
\mathcal{K}(S_{n-1}) = n-1$ et que $\mathcal{K}(S_{n-1}) \subset
\mathcal{K}(S)$. Le lemme est alors imm\'{e}diat. \qed

\medskip

{\bf 5.2.}~{\bf Th\'{e}or\`{e}me.}~{\em Soit $\mathfrak{g}$ une
alg\`{e}bre semi-simple de rang $\ell$. Pour tout entier $i$
v\'{e}rifiant $0 \leqslant i \leqslant \ell$, il existe une
sous-alg\`{e}bre parabolique $\mathfrak{p}$ de $\mathfrak{g}$ telle
que $\chi (\mathfrak{p}) = i$.}

\smallskip

{\em Preuve.}~Il est clair que l'on peut supposer $\mathfrak{g}$
simple. Soit $n = \card \mathcal{K}(\Pi )$. On pose $\Pi
=\{\alpha_{1}, \dots , \alpha_{\ell}\}$ en utilisant la
num\'{e}rotation des syst\`{e}mes de racines de [1], pages 250 \`{a}
275.

Il existe des parties $S_{0}, S_{1}, \dots , S_{n}$ de $\Pi$, o\`{u}
$S_{0} = \emptyset$ et $S_{n} = \Pi$, telles que $S_{1}, \dots ,
S_{n}$ v\'{e}rifient les conditions du lemme 5.1.

D'apr\`{e}s la proposition 4.4, pour $0 \leqslant i \leqslant n$, on a
:
$$
\chi (\mathfrak{g}_{\Pi, S_{i}}) = \ell + i - n.
$$

Le th\'{e}or\`{e}me est donc \'{e}tabli si $\mathfrak{g}$ est de l'un
des types $B_{k}$, $C_{k}$, $D_{2k}$, $E_{7}$, $E_{8}$, $F_{4}$ ou
$G_{2}$ car, dans ces cas, $n = \ell$ d'apr\`{e}s le tableau 1.6.

1) Supposons $\mathfrak{g}$ de type $D_{2k+1}$. Alors $\card
\mathcal{K}(\Pi ) = 2k$ et, si $1 \leqslant i \leqslant \ell$, ce qui
pr\'{e}c\`{e}de montre qu'il existe une sous-alg\`{e}bre parabolique
de $\mathfrak{g}$ dont l'indice est \'{e}gal \`{a} $i$. D'autre part,
on voit facilement que $\alpha_{2k} \notin E_{\Pi}$. Compte tenu de
l'assertion (iii) de 4.4, il vient $\chi (\mathfrak{g}_{\Pi ,\{
\alpha_{2k}\}}) = 0$.

2) Si $\mathfrak{g}$ est de type $E_{6}$, on a $\card \mathcal{K}(\Pi
) = 4$. Pour $2 \leqslant i \leqslant 6$, il existe donc une
sous-alg\`{e}bre parabolique de $\mathfrak{g}$ dont l'indice est
\'{e}gal \`{a} $i$.

On v\'{e}rifie ais\'{e}ment que $\{ \alpha_{1}, \alpha_{5}\} \cup \{
\varepsilon_{K}\, ; \, K \in \mathcal{K}(\Pi )\}$ est une base de
$\mathfrak{h}^{*}$. A nouveau d'apr\`{e}s 4.4, on obtient :
$$
\chi (\mathfrak{g}_{\Pi , \{ \alpha_{1}\}}) = 1 \ , \ \chi
(\mathfrak{g}_{\Pi , \{ \alpha_{1}, \alpha_{5}\}}) = 0.
$$

3) Supposons $\mathfrak{g}$ de type $A_{\ell}$, et soit $\rho$ la
partie enti\`{e}re de $\dfrac{\ell +1}{2}$. D'apr\`{e}s 1.6, $\rho$
est le cardinal de $\mathcal{K}(\Pi )$. Si $\ell' = \ell - \rho$, il
nous faut prouver l'existence d'une sous-alg\`{e}bre parabolique de
$\mathfrak{g}$ d'indice $0, 1, \dots , \ell' -1$.

Si $\ell = 2t$ est pair, on a 
$$
\mathcal{K}(\Pi ) = \{ \{ \alpha_{j+1}, \dots , \alpha_{2t-j}\} \, ;
\, 0 \leqslant j \leqslant t-1\} 
$$
et, si $\ell = 2t+1$ est impair, alors :
$$
\mathcal{K}(\Pi ) = \{ \{ \alpha_{j+1}, \dots , \alpha_{2t+1-j}\} \, ;
\, 0 \leqslant j \leqslant t\} .
$$

Posons $T_{0} = \emptyset$, et d\'{e}finissons $T_{k}$, pour $1
\leqslant k \leqslant \ell'$, par :
$$
T_{k} = \begin{cases}
T_{k-1} \cup \{ \alpha_{k}\} \ \text{ si } \ k \ \text{ est impair,}
\hfill \\
T_{k-1} \cup \{ \alpha_{\ell + 1 -k}\} \ \text{ si } \ k \ \text{ est
pair.}
\end{cases}
$$

On v\'{e}rifie facilement que les $\varepsilon_{M}$, pour $M \in
\mathcal{K}(\Pi ) \cup \mathcal{K}(T_{k})$ forment une famille
libre. Compte tenu de 4.4\, (iii), il vient $\chi
(\mathfrak{g}_{\Pi , T_{k}}) = \ell' - k$ si $0 \leqslant k
\leqslant \ell'$. Comme $\mathfrak{g}_{\Pi , T_{k}}$ est une
sous-alg\`{e}bre parabolique de $\mathfrak{g}$, on a obtenu le
r\'{e}sultat. \qed

\medskip

{\bf Remarque.}~Le th\'{e}or\`{e}me 5.2 est \'{e}nonc\'{e} sans
d\'{e}monstration dans [4].

\vspace{5mm}

{\bf Bibliographie}

\medskip

[1]~{\sc Bourbaki N.}, {\em Groupes et alg\`{e}bres de Lie},
chap. 4,5,6, Masson, 1981.

[2]~{\sc Dergachev V., Kirillov A.}, Index of Lie algebras of seaweed
type, {\em J. of Lie Theory}, 10, 2000, p. 331-343

[3]~{\sc Dixmier J.}, {\em Enveloping algebras}, Graduate Studies in
Math., 11, AMS, 1996.

[4]~{\sc Elashvili A. G.}, On the index of parabolic subalgebras of
semisimple Lie algebras, {\em Preprint}, 1990.

[5]~{\sc Jantzen J. C.}, {\em Einh\"{u}llenden Algebren halbeinfacher
Lie-Algebren}, Ergebnisse der Mathematik and iher Grenzgebiete, 3,
Springer-Verlag, 1983.

[6]~{\sc Joseph A.}, A preparation theorem for the prime spectrum of a
semisimple Lie algebra, {\em J. of Algebra}, 48, 1977, p. 241-289.

[7]~{\sc Panysushev D.,} Inductive formulas for the index of seaweed
Lie algebras, {\em Moscow Math. Journal}, 1, 2001, p. 221-241.

[8]~{\sc Tauvel P.}, {\em Introduction \`{a} la th\'{e}orie des
alg\`{e}bres de Lie}, Diderot, 1998.

[9]~{\sc Tauvel P., Yu R.W.T.}, Indice et formes lin\'{e}aires stables
dans les alg\`{e}bres de Lie, {\em \`{a} para\^{\i}tre dans J. of
Algebra}. 

\bigskip

Patrice Tauvel et Rupert W.T. Yu

Universit\'{e} de Poitiers

UMR 6086 du CNRS, D\'{e}partement de math\'{e}matiques

T\'{e}l\'{e}port 2 - BP 30179

Boulevard Marie et Pierre Curie

86962 Futuroscope Chasseneuil Cedex France

tauvel@math.univ-poitiers.fr 

yuyu@math.univ-poitiers.fr

\end{document}